\documentclass[12pt]{article}

\usepackage{latexsym,amsmath}
\usepackage{amssymb}

\newtheorem{thm}{Theorem}[section]

\newtheorem{lemma}[thm]{Lemma}
\newtheorem{prop}[thm]{Proposition}
\newtheorem{cor}[thm]{Corollary}

\newcommand{\beq}[1]{\begin{equation}\label{#1}}
\newcommand{\enq}[0]{\end{equation}}

\newcommand{\qed}[0]{{\hspace*{\fill}\mbox{$\Box$}}}

\begin{document}
\renewcommand{\thefootnote}{\fnsymbol{footnote}}

\title{On weighted graph homomorphisms}

\author{David Galvin\thanks{Microsoft Research, 1 Microsoft Way, Redmond WA 98052, galvin@microsoft.com.}\\
Prasad Tetali\thanks{School of Mathematics \& College of
Computing, Georgia Institute of Technology, Atlanta, GA
30332-0160. Research supported in part by NSF grant DMS-0100289.}}

\date{Appeared 2004}

\maketitle

\begin{abstract}
For given graphs $G$ and $H$, let $|Hom(G,H)|$ denote the set of
graph homomorphisms from $G$ to $H$. We show that for any finite,
$n$-regular, bipartite graph $G$ and any finite graph $H$ (perhaps
with loops), $|Hom(G,H)|$ is maximum when $G$ is a disjoint union
of $K_{n,n}$'s. This generalizes a result of J. Kahn on the number
of independent sets in a regular bipartite graph. We also give the
asymptotics of the logarithm of $|Hom(G,H)|$ in terms of a simply
expressed parameter of $H$.

We also consider weighted versions of these results which may be
viewed as statements about the partition functions of certain
models of physical systems with hard constraints.
\end{abstract}

\section{Introduction}

Let $G$ be an $n$-regular, $N$-vertex bipartite graph on vertex
set $V(G)$, and let $H$ be a fixed graph on vertex set $V(H)$
(perhaps with loops). We will always use $u$, $v$ for the vertices
of $G$ and $i$,$j$ for those of $H$. Set
$$Hom(G,H)=\{f:V(G)\rightarrow V(H)~:~u\sim v\Rightarrow f(u)\sim
f(v)\}.$$ That is, $Hom(G,H)$ is the set of graph homomorphisms
from $G$ to $H$. (For graph theory basics, see e.g.
\cite{MGT},~\cite{Diestel}).

When $H=H_{ind}$ consists of one looped and one unlooped vertex
connected by an edge, an element of $Hom(G,H_{ind})$ can be
thought of as a specification of an independent set (a set of
vertices spanning no edges) in $G$. Our point of departure is the
following result of Kahn \cite{JK}, bounding the number of
independent sets in regular bipartite graphs. For any graph $G$,
write ${\cal I}(G)$ for the set of independent sets of $G$.

\begin{thm} \label{jeff}
For any $n$-regular, $N$-vertex bipartite graph $G$,
$$|{\cal I}(G)| \leq (2^{n+1}-1)^{N/2n}.$$
\end{thm}

An approximate version of Theorem \ref{jeff} --- $\log |{\cal
I}(G)| \leq (1/2 +o(1))N$, where $o(1) \rightarrow 0$ as $n
\rightarrow \infty$ --- for general $n$-regular, $N$-vertex $G$
was earlier proved by Alon \cite{Noga}. Note that
$|Hom(K_{n,n},H_{ind})|=2^{n+1}-1$ (where $K_{n,n}$ is the
complete bipartite graph with $n$ vertices on each side), so we
may paraphrase Theorem \ref{jeff} by saying that
$|Hom(G,H_{ind})|$ is maximum when $G$ is a disjoint union of
$K_{n,n}$'s. Our main result is a generalization of this statement
(and our proof is a generalization of Kahn's).

\begin{prop} \label{ub}
For any $n$-regular, $N$-vertex bipartite $G$, and any $H$,
$$|Hom(G,H)|\leq |Hom(K_{n,n},H)|^{N/2n}.$$
\end{prop}

Somewhat surprisingly, we can also exhibit a lower bound that is
good enough to allow us to obtain the asymptotics of $\log
|Hom(G,H)|$ for fixed $H$ as $n \rightarrow \infty$ (here, and
throughout the rest of the paper, we use $\log$ for the base $2$
logarithm). To state the result, it is convenient to introduce a
parameter of $H$ that is very closely related to
$|Hom(K_{n,n},H)|$, but is easier to work with. Set
$$\eta(H)=\max\{|A||B|:A,B\subseteq V(H),i\sim j~\forall i\in A,j\in B\}.$$
(When $H$ is loopless, this is the maximum number of edges in a
complete bipartite subgraph of $H$. Peeters \cite{peeters} has
recently shown that determining $\eta(H)$, even when $H$ is
bipartite, is NP-complete.)

\begin{prop} \label{lb}
For any $n$-regular, $N$-vertex bipartite $G$, and any $H$,
$$\frac{\log\eta(H)}{2} \leq \frac{\log|Hom(G,H)|}{N} \leq \frac{\log\eta(H)}{2} + \frac{|V(H)|}{2n}.$$
\end{prop}

We use the example of $H=K_k$, the complete graph on $k$ vertices,
to illustrate the definition of $\eta$. It is easy to see that for
any $A,B \subseteq V(K_k)$, we have $i\sim j ~\forall i\in A, j\in
B$ iff $A$ and $B$ are disjoint, and so $|A||B|$ is maximum when
$|A|$ and $|B|$ are as close as possible to $k/2$. Hence
$\eta(K_k)=\lfloor k/2 \rfloor \lceil k/2 \rceil$. Since an
element of $Hom(G,K_k)$ is exactly a proper $k$ coloring of $G$,
we get as a corollary of Proposition \ref{lb} an approximate count
of the number of $k$-colorings of a regular bipartite graph.

\begin{cor} \label{col}
For any $n$-regular, $N$-vertex bipartite $G$,
$$|Hom(G,K_k)| = \left(\lfloor k/2 \rfloor \lceil k/2 \rceil\right)^{N\left(1/2+o(1)\right)}.$$
\end{cor}

\medskip

We now consider weighted versions of Propositions \ref{ub} and
\ref{lb}. Following \cite{BW}, we put a measure on $Hom(G,H)$ as
follows. To each $i\in V(H)$ assign a positive ``activity''
$\lambda_i$, and write $\Lambda$ for the set of activities. Give
each $f \in Hom(G,H)$ weight
$$w^{\Lambda}(f)=\prod_{v\in V(G)}\lambda_{f(v)}.$$
The constant that turns this assignment of weights on $Hom(G,H)$
into a probability distribution is
$$Z^{\Lambda}(G,H)=\sum_{f\in Hom(G,H)}w^{\Lambda}(f).$$
When all activities are $1$, we have
$Z^{\Lambda}(G,H)=|Hom(G,H)|$, and so the following is a
generalization of Proposition \ref{ub}.

\begin{prop} \label{partition2}
For any $n$-regular, $N$-vertex bipartite $G$, any $H$, and any
system $\Lambda$ of positive activities on $V(H)$,
$$Z^{\Lambda}(G,H) \leq \left(Z^{\Lambda}(K_{n,n},H)\right)^{N/2n}.$$
\end{prop}

It was observed in \cite{BW} that $Z^{\Lambda}(G,H)$ may be
related to $|Hom(G,H')|$ for an appropriate modification $H'$ of
$H$. That observation (which will be discussed in more detail in
Section \ref{sectionproofs}) is central to the proof of
Proposition \ref{partition2}.

Proposition \ref{lb} also generalizes. For a set of activities
$\Lambda$ on $V(H)$, set
$$\eta^{\Lambda}(H)=\max\left\{\left(\sum_{i\in
A}\lambda_i\right)\left(\sum_{j\in B}\lambda_j\right):A,B
\subseteq V(H), i\sim j~\forall i\in A, j\in B\right\}.$$

\begin{prop} \label{partition}
For any $n$-regular, $N$-vertex bipartite $G$, any $H$, and any
system $\Lambda$ of positive activities on $V(H)$,
$$\frac{\log\eta^{\Lambda}(H)}{2} \leq \frac{\log
Z^{\Lambda}(G,H)}{N} \leq \frac{\log\eta^{\Lambda}(H)}{2} +
\frac{|V(H)|}{2n}.$$
\end{prop}

\medskip

We may put these results in the framework of a well-known
mathematical model of physical systems with ``hard constraints"
(see \cite{BW}). These are systems with strictly forbidden
configurations. An example is the hard-core lattice gas model, in
which a legal configuration of particles on a lattice is precisely
one in which no two adjacent lattice sites are occupied. (By way
of contrast, consider the ferromagnetic Ising model, where
adjacent particles are discouraged from having opposing spins, but
not forbidden --- this is a ``soft constraint".)

We think of the vertices of $G$ as particles and the edges as
bonds between pairs of particles, and we think of the vertices of
$H$ as possible ``spins'' that particles may take. Pairs of
vertices of $G$ joined by a bond may have spins $i$ and $j$ only
when $i$ and $j$ are adjacent in $H$ (in particular, they may both
have spin $i$ only when $i$ has a loop in $H$). Thus the legal
spin configurations on the vertices of $G$ are precisely the
homomorphisms from $G$ to $H$. We think of the activities on the
vertices of $H$ as a measure of the likelihood of seeing the
different spins; the probability of a particular spin
configuration is proportional to the product over the vertices of
$G$ of the activities of the spins. Propositions \ref{partition2}
and \ref{partition} concern the ``partition function'' of this
model --- the normalizing constant that turns the above-described
system of weights on the set of legal configurations into a
probability measure.

\medskip

The results we actually prove are in a slightly more general
weighted model. Write ${\cal E}_G$ and ${\cal O}_G$ for the
partition classes of $G$, and to each $i\in V(H)$ assign a
positive {\em pair} of activities $(\lambda_i,\mu_i)$. Write
$(\Lambda,\rm{M})$ for the set of activities. Give each $f \in
Hom(G,H)$ weight
$$w^{(\Lambda,\rm{M})}(f)=\prod_{v\in {\cal E}_G}\lambda_{f(v)}
\prod_{v\in {\cal O}_G}\mu_{f(v)}.$$ The constant that turns this
assignment of weights on $Hom(G,H)$ into a probability
distribution is \beq{lm}Z^{(\Lambda,\rm{M})}(G,H)=\sum_{f\in
Hom(G,H)}w^{(\Lambda,\rm{M})}(f).\enq

A special case of this model was considered by Kahn \cite{JK2}
(see also \cite{Haggstrom}), where Theorem \ref{jeff} was extended
to

\begin{thm} \label{jeff2}
For any $n$-regular, $N$-vertex bipartite $G$, and any $\lambda,
\mu \geq 1$,
$$\sum_{I \in {\cal I}(G)} \prod_{v\in {\cal E}_G} \lambda^{|I\cap {\cal E}_G|}
\prod_{v\in {\cal O}_G} \mu^{|I\cap {\cal O}_G|} \leq
((1+\lambda)^n+(1+\mu)^n-1)^{N/2n}.$$
\end{thm}

It was conjectured in \cite{JK2} that the assumption $\lambda, \mu
\geq 1$ may be relaxed to $\lambda, \mu \geq 0$. We show that this
is indeed true, by generalizing Proposition \ref{partition2} to:

\begin{prop} \label{partition3}
For any $n$-regular, $N$-vertex bipartite $G$, any $H$, and any
system $(\Lambda,\rm{M})$ of positive activities on $V(H)$,
$$Z^{(\Lambda,\rm{M})}(G,H) \leq \left(Z^{(\Lambda,\rm{M})}(K_{n,n},H)\right)^{N/2n}.$$
\end{prop}

We also generalize Proposition \ref{partition} to this setting.
Set
$$\eta^{(\Lambda,\rm{M})}(H)=
\max\left\{\left(\sum_{i\in A}\lambda_i\right)\left(\sum_{j\in
B}\mu_j\right):A,B \subseteq V(H), i\sim j~\forall i\in A, j\in
B\right\}.$$

\begin{prop} \label{partition4}
For any $n$-regular, $N$-vertex bipartite $G$, any $H$, and any
system $(\Lambda,\rm{M})$ of positive activities on $V(H)$,
$$\frac{\log\eta^{(\Lambda,\rm{M})}(H)}{2} \leq \frac{\log Z^{(\Lambda,\rm{M})}(G,H)}{N}
\leq \frac{\log\eta^{(\Lambda,\rm{M})}(H)}{2} +
\frac{|V(H)|}{2n}.$$
\end{prop}

Proposition \ref{partition3} generalizes to the case of biregular
$G$ (a bipartite graph $G$ with partition classes ${\cal E}_G$ and
${\cal O}_G$ is {\em $(a,b)$-biregular} if all vertices in ${\cal
E}_G$ have degree $a$ and all in ${\cal O}_G$ have degree $b$).
The proof of the following proposition, which is a straightforward
modification of the proof of Proposition \ref{partition3}, is
omitted.

\begin{prop} \label{bireg}
For any $(a,b)$-biregular, $N$-vertex, bipartite $G$, any $H$, and
any system $(\Lambda,\rm{M})$ of positive activities on $V(H)$,
$$Z^{(\Lambda,\rm{M})}(G,H) \leq \left(Z^{(\Lambda,\rm{M})}(K_{a,b},H)\right)^{N/(a+b)}.$$
\end{prop}

It was conjectured in \cite{JK} that Theorem \ref{jeff} remains
true without the assumption that $G$ is bipartite. We similarly
conjecture that biparticity is unnecessary in Proposition
\ref{partition3}, and hence also in Propositions \ref{ub} and
\ref{partition2}\footnote{Note added for ArXiv submission: Propositions \ref{ub} turns out not to be true for general $H$ without the assumption that $G$ is bipartite. See D. Galvin, Maximizing $H$-colorings of regular graphs, {\em J. Graph Theory} \& arXiv:1110.3758, for a discussion of an amended conjecture}. (Proposition \ref{lb}, and hence also
Propositions \ref{partition} and \ref{partition4}, is easily seen
to fail for non-bipartite $G$.)

\medskip

The proof of Proposition \ref{partition3} requires entropy
considerations; these are reviewed in Section
\ref{sectionentropy}. The proofs are then given in Section
\ref{sectionproofs}.

\section{Entropy} \label{sectionentropy}
\label{Ent}

Here we briefly review the relevant entropy material. Our
treatment is mostly copied from \cite{JK}. For a more thorough
discussion, see e.g. \cite{McE}.

In what follows ${\bf X}$, ${\bf Y}$ etc. are discrete random
variables, which in our usage are allowed to take values in any
finite set.

The {\em entropy} of the random variable ${\bf X}$ is
$$H({\bf X}) = \sum_x p(x)\log\frac{1}{p(x)},$$
where we write $p(x)$ for ${\bf P}({\bf X}=x)$ (and extend this
convention in natural ways below). The {\em conditional entropy}
of ${\bf X}$ given ${\bf Y}$ is
$$H({\bf X}|{\bf Y}) ={\bf E} H({\bf X}|\{{\bf Y}=y\})
=\sum_yp(y)\sum_xp(x|y)\log\frac{1}{p(x|y)}.
$$
Notice that we are also writing $H({\bf X}|Q)$ with $Q$ an event
(in this case $Q=\{{\bf Y} =y\}$):
$$H({\bf X}|Q)=\sum p(x|Q)\log\frac{1}{p(x|Q)}.$$
When we condition on a random variable and an event
simultaneously, we use ``;'' to separate the two.

For a random vector ${\bf X}=({\bf X}_1,\ldots ,{\bf X}_n)$ (note
this is also a random variable), we have
\begin{equation} \label{chain}
H({\bf X}) =H({\bf X}_1)+H({\bf X}_2|{\bf X}_1)+\cdots
    +H({\bf X}_n|{\bf X}_1,\ldots, {\bf X}_{n-1}).
\end{equation}

We will make repeated use of the inequalities
\begin{equation} \label{range}
H({\bf X}) \leq \log |{\rm range}({\bf X})|~~~~ \mbox{(with
equality if ${\bf X}$ is uniform),}
\end{equation}
$$H({\bf X}|{\bf Y}) \leq H({\bf X}),$$ and more generally,
\begin{equation} \label{dropping}
\mbox{if ${\bf Y}$ determines ${\bf Z}$ then $H({\bf X}|{\bf Y})
\leq H({\bf X}|{\bf Z})$.}
\end{equation}

Note that (\ref{chain}) and (\ref{dropping}) imply
$$H({\bf X}) \leq H({\bf Y}) + H({\bf X}|{\bf Y})$$ and
\begin{equation} \label{sub}
H({\bf X}_1,\ldots, {\bf X}_n)\leq \sum H({\bf X}_i)
\end{equation}
We also have a conditional version of (\ref{sub}):
$$H({\bf X}_1,\ldots, {\bf X}_n|{\bf Y})\leq \sum H({\bf X}_i|{\bf Y}).$$

We will also need the following lemma of Shearer (see \cite[p.
33]{CFGS}). For a random vector ${\bf X}=({\bf X}_1,\ldots, {\bf
X}_m)$ and $A\subseteq [m]$, set ${\bf X}_A=({\bf X}_i:i\in A)$.

\begin{lemma}
\label{Lshearer} Let ${\bf X}=({\bf X}_1,\ldots, {\bf X}_m)$ be a
random vector and ${\cal A}$ a collection of subsets (possibly
with repeats) of $[m]$, with each element of $[m]$ contained in at
least $t$ members of ${\cal A}$. Then
$$H({\bf X})\leq
    \frac{1}{t}\sum_{A\in{\cal A}}H({\bf X}_A).$$
\end{lemma}

\section{Proofs} \label{sectionproofs}

We begin by setting up some conventions. For a regular, bipartite
graph $G$, we write ${\cal E}_G$ and ${\cal O}_G$ for the
partition classes. For ease of notation, we write ${\cal E}_n$ for
${\cal E}_{K_{n,n}}$ and ${\cal O}_n$ for ${\cal O}_{K_{n,n}}$.

For a partition $U \cup L$ of $V(H)$, set
$$Hom^{U,L}(G,H)=\{f \in Hom(G,H):f({\cal E}_G) \subseteq U, f({\cal O}_G) \subseteq
L\}.$$ (For a set $X$ we write $f(X)$ for $\{f(x):x\in X\}$.)

We begin by deriving a useful expression for
$|Hom^{U,L}(K_{n,n},H)|$. For each $A \subseteq L$ set
$${\cal H}(A) = \{f \in Hom^{U,L}(K_{n,n},H):f({\cal O}_n)=A\},$$
$$T(A) = \{g:[n]\rightarrow A~:~\mbox{$g$ surjective}\}$$
and
$$C^U(A)=\{j \in U:j\sim i~\forall i\in A\}.$$
(Observe that $C^U(A)$ is the set of all possible images of $v \in
{\cal E}_G$ under a member of $Hom^{U,L}(G,H)$, given that the
image of $N(v)$ is $A$.) It is easy to see that $\{{\cal H}(A):A
\subseteq L\}$ forms a partition of $Hom^{U,L}(K_{n,n},H)$, and
also that for each $A$, $|{\cal H}(A)|=|T(A)||C^U(A)|^n$. Thus we
have
\begin{equation} \label{knn'}
|Hom^{U,L}(K_{n,n},H)| = \sum_{A\subseteq L} |T(A)||C^U(A)|^n.
\end{equation}

The following is the central lemma in the proofs of Propositions
\ref{partition3} and \ref{partition4}. The proof is based on
\cite[Thm. 1.9]{JK}.

\begin{lemma} \label{ub'}
For any $n$-regular, $N$-vertex bipartite $G$, and any $H$ with
$U\cup L$ a partition of $V(H)$,
$$|Hom^{U,L}(G,H)|\leq |Hom^{U,L}(K_{n,n},H)|^{N/2n}.$$
\end{lemma}

\noindent {\em Proof: }Let ${\bf f}$ be chosen uniformly from
$Hom^{U,L}(G,H)$. For $v \in V(G)$, write ${\bf f}_v$ for ${\bf
f}(v)$, ${\bf N}_v$ for ${\bf f}|_{N(v)}$ and ${\bf M}_v$ for
$\{{\bf f}_w:w\in N(v)\}$. For $v \in {\cal E}_G$ and $A\subseteq
L$, write $m_v(A)$ for ${\bf P}({\bf M}_v=A)$. (Note that $\sum_A
m_v(A) = 1$.) We have (with the main inequalities justified below;
the remaining steps follow in a straightforward way from the
material of Section \ref{sectionentropy})
\begin{eqnarray}
\log|Hom^{U,L}(G,H)| & = & H({\bf f}) \nonumber \\
& = & H({\bf f}|_{{\cal O}_G})+H({\bf f}|_{{\cal E}_G}~|~{\bf
f}|_{{\cal O}_G})
\nonumber \\
& \leq & H({\bf f}|_{{\cal O}_G}) + \sum_{v\in {\cal E}_G} H({\bf
f}_v~|~{\bf f}|_{{\cal O}_G})
\nonumber \\
& \leq & H({\bf f}|_{{\cal O}_G}) + \sum_{v\in {\cal E}_G} H({\bf
f}_v~|~{\bf N}_v)
\nonumber \\
& \leq & \frac{1}{n}\sum_{v\in {\cal E}_G} H({\bf N}_v) +
\sum_{v\in {\cal E}_G} H({\bf f}_v~|~{\bf N}_v) \label{shearer}\\
& \leq & \frac{1}{n}\sum_{v\in {\cal E}_G} \left[H({\bf
M}_v)+H({\bf N}_v|{\bf M}_v)\right] + \sum_{v\in {\cal E}_G}
H({\bf
f}_v~|~{\bf N}_v) \nonumber \\
& \leq & \frac{1}{n}\sum_{v\in {\cal E}_G}
\left[H({\bf M}_v)+H({\bf N}_v|{\bf M}_v)+nH({\bf f}_v|{\bf N}_v)\right] \nonumber \\
& \leq & \frac{1}{n}\sum_{v\in {\cal E}_G} \sum_{A\subseteq L}
\Big[m_v(A)\log\frac{1}{m_v(A)} + \nonumber \\
& & ~~~~~~~~~~~~~ m_v(A)H({\bf N}_v|\{{\bf M}_v=A\})+  \nonumber \\
& & ~~~~~~~~~~~~~~~~~ nm_v(A)H({\bf f}_v|{\bf N}_v;\{{\bf
M}_v=A\})\Big] \nonumber
\end{eqnarray}
\begin{eqnarray}
~~~~~~~~~~~~~~~~~~~~~ & \leq & \frac{1}{n}\sum_{v\in {\cal E}_G}
\sum_{A\subseteq L}
\Big[m_v(A)\log\frac{1}{m_v(A)} + \nonumber \\
& & ~~~~~~~ m_v(A)\log|T(A)| + \nonumber \\
& & ~~~~~~~~~~~~nm_v(A)\log|C^U(A)|\Big] \label{usingrange} \\
& = & \frac{1}{n}\sum_{v\in {\cal E}_G} \sum_{A\subseteq L}
m_v(A)\log\frac{|T(A)||C^U(A)|^n}{m_v(A)} \nonumber \\
& \leq & \frac{1}{n}\sum_{v\in {\cal E}_G}
\log \left[ \sum_{A\subseteq L} |T(A)||C^U(A)|^n \right] \label{jensen} \\
& = & \frac{N}{2n} \log|Hom^{U,L}(K_{n,n},H)|. \label{usingknn'}
\end{eqnarray}
The main inequality (\ref{shearer}) involves an application of
Lemma \ref{Lshearer}, with ${\cal A}=\{N(v):v \in {\cal E}_G\}$,
and (\ref{jensen}) is an application of Jensen's inequality. In
(\ref{usingrange}), we use (\ref{range}), noting that conditioning
on the event $\{{\bf M}_v = A\}$ there are $|T(A)|$ possible
values for ${\bf N}_v$, and $|C^U(A)|$ possible values for ${\bf
f}_v$. Finally, (\ref{usingknn'}) follows from (\ref{knn'}). \qed

It is worth noting at this point that Lemma \ref{ub'} easily
implies Proposition \ref{ub}. Let $H'$ be the graph on vertex set
$\cup_{i \in V(H)} \{v_i,w_i\}$ with $v_i$ and $w_j$ adjacent
exactly when $i$ and $j$ are adjacent in $H$. Set $U=\{v_1,
\ldots, v_{|V(H)|}\}$ and $L=\{w_1, \ldots, w_{|V(H)|}\}$. It is
easy to check that $|Hom(G,H)|=|Hom^{U,L}(G,H')|$, from which
Proposition \ref{ub} follows via an application of Lemma
\ref{ub'}.

This idea of ``doubling'' $H$, combined with the construction of
\cite{BW} that relates $Z^{\Lambda}(G,H)$ to $|Hom(G,H')|$ for an
appropriate modification $H'$ of $H$, allows us to pass from
Proposition \ref{ub} to Proposition \ref{partition3}. The details
are as follows.

Recall that our aim is to upper bound the partition function
$Z^{(\Lambda,\rm{M})}(G,H)$ (see (\ref{lm})). By continuity, we
may assume that all activities are rational. Let $C$ be the least
positive integer such that $C\lambda_i$ and $C\mu_i$ are integers
for each $i\in V(H)$. Let $H^{(\Lambda,\rm{M})}$ be the graph
whose vertex set is obtained from $H$ by replacing each $i\in
V(H)$ by two sets, $D^U_i = \{i^U_1, \ldots ,i^U_{C\lambda_i}\}$
and $D^L_i = \{i^L_1, \ldots ,i^L_{C\mu_i}\}$ of $C\lambda_i$ and
$C\mu_i$ vertices. For each $i, j\in V(H)$ (not necessarily
distinct), $i' \in D^U_i$ and $j' \in D^L_j$, join $i'$ to $j'$
exactly when $i$ and $j$ are adjacent in $H$. Set
$U=U(H^{(\Lambda,\rm{M})})=\cup_{i \in V(H)}D^U_i$ and
$L=L(H^{(\Lambda,\rm{M})})=\cup_{i \in V(H)}D^L_i$.

We wish to relate $Z^{(\Lambda,\rm{M})}(G,H)$ to
$|Hom^{U,L}(G,H^{(\Lambda,\rm{M})})|$. Say that a function $g \in
Hom^{U,L}(G,H^{(\Lambda,\rm{M})})$ is a {\em lift} of $f\in
Hom(G,H)$ if for all $v \in V(G)$,
$$g(v)=
\left\{\begin{array}{ll}
               f(v)^U_k,\mbox{ some } 1\leq k\leq
C\lambda_{f(v)} & \mbox{if $v\in {\cal E}_G$,}\\
               f(v)^L_k,\mbox{ some } 1\leq k\leq C\mu_{f(v)} & \mbox{if
$v\in {\cal O}_G$.}
\end{array}\right.$$
Set
$${\cal G}(f)=\{g \in Hom^{U,L}(G,H^{(\Lambda,\rm{M})}):\mbox{$g$ is a
lift of $f$}\}.$$ It is easy to check that $|{\cal
G}(f)|=w^{(\Lambda,\rm{M})}(f)C^N$ for each $f \in Hom(G,H)$, and
that $\{{\cal G}(f):f \in Hom(G,H)\}$ forms a partition of
$Hom^{U,L}(G,H^{(\Lambda,\rm{M})})$. From this it follows that
\begin{equation} \label{unweighting'}
Z^{(\Lambda,\rm{M})}(G,H)=\frac{|Hom^{U,L}(G,H^{(\Lambda,\rm{M})})|}{C^N}.
\end{equation}

\medskip

We now have all we need to prove Propositions \ref{partition3} and
\ref{partition4}.

\noindent {\em Proof of Proposition \ref{partition3}: }Applying
(\ref{unweighting'}) with $G=K_{n,n}$ we get
\begin{equation} \label{un2}
Z^{(\Lambda,\rm{M})}(K_{n,n},H) =
\frac{|Hom^{U,L}(K_{n,n},H^{(\Lambda,\rm{M})})|}{C^{2n}}.
\end{equation}
Proposition \ref{partition3} now follows from
(\ref{unweighting'}), (\ref{un2}) and Lemma \ref{ub'}.\qed

\medskip

\noindent {\em Proof of Proposition \ref{partition4}: }For each $A
\subseteq V(H)$, set $C(A)=\{j\in H:j \sim i~\forall i\in A\}$ and
$${\cal D}(A)=\{f\in Hom(K_{n,n},H):f({\cal E}_n)\subseteq A,
f({\cal O}_n)\subseteq C(A)\}.$$ By Proposition \ref{partition3}
we have
\begin{eqnarray}
\left(Z^{(\Lambda,\rm{M})}(G,H)\right)^{2n/N} & \leq &
Z^{(\Lambda,\rm{M})}(K_{n,n},H) \nonumber \\
& \leq & \sum_{A\subseteq V(H)} \sum_{f\in {\cal D}(A)}
w^{(\Lambda,\rm{M})}(f) \nonumber \\
& = & \sum_{A\subseteq V(H)} \left(\sum_{i\in
A}\lambda_i\right)^n\left(\sum_{j \in C(A)}\mu_j\right)^n
\nonumber \\
& \leq & 2^{|V(H)|}\left(\eta^{(\Lambda,\rm{M})}(H)\right)^n.
\end{eqnarray}
This gives the upper bound. For the lower bound, let $A,B
\subseteq V(H)$ satisfying $i\sim j~\forall i\in A, j\in B$ be
such that
$$\eta^{(\Lambda,\rm{M})}(H)=\left(\sum_{i\in
A}\lambda_i\right)\left(\sum_{j\in B}\mu_j\right).$$ We have
\begin{eqnarray*}
Z^{(\Lambda,\rm{M})}(G,H) & \geq & \sum\{w^{(\Lambda,\rm{M})}(f)
:f({\cal E}_G)\subseteq
A, f({\cal O}_G)\subseteq B\} \\
& = & \left(\sum_{i\in
A}\lambda_i\right)^{N/2}\left(\sum_{j\in B}\mu_j\right)^{N/2} \\
& = & \eta^{(\Lambda,\rm{M})}(H)^{N/2}.
\end{eqnarray*}
\qed

\bigskip \noindent
{\bf Acknowledgment} This research was done while the second
author was visiting Microsoft Research. He would like to thank
Microsoft Research, and especially the theory group, for providing
him with this opportunity.

\end{document}